\documentclass[10pt, a4paper]{article}

\IfFileExists{marvosym.sty}{
    \usepackage{marvosym}
    \newcommand{\emailsymbol}{\Letter}
}{
    \newcommand{\emailsymbol}{\dagger}
}
\makeatletter
\def\@maketitle{%
  \newpage
  \null
  \vskip 1.5em
  \begin{flushleft}
    {\LARGE\bfseries \@title \par}
    \vskip 1em
    {\large \@author \par}
    \vskip 0.6em
    {\small \@date \par}
  \end{flushleft}
  \vskip 0.8em
}
\makeatother
\title{\textsc{A Remark on the Topology of the Regular Loci of Some Complexified Hamiltonian Systems}}
\author{\textbf{\textsc{Zhiyuan Liu}}\textsuperscript{1,2,\emailsymbol}}
\date{}
\usepackage[
  left=2.5cm,      
  right=2.5cm,     
  top=2cm,       
  bottom=2cm,    
  footskip=1cm     
]{geometry}
\usepackage{placeins}
\usepackage{amsfonts}
\usepackage{amssymb}
\usepackage{amsmath,amscd}
\usepackage{bm}
\usepackage{slashed}
\usepackage{multirow}
\usepackage{makecell}
\usepackage{mathrsfs} 
\usepackage{tikz-cd}  
\usetikzlibrary{decorations.pathreplacing}
\usetikzlibrary{3d, calc, fadings}
\usepackage{extarrows} 
\newcommand{\R}{{\mathbb R}}
\newcommand{\C}{{\mathbb C}}
\newcommand{\Z}{{\mathbb Z}}
\newcommand{\h}{{\mathbb H}}

\usepackage{indentfirst} 
\usepackage{sectsty}
\usepackage{fancyhdr} 
\usepackage{graphicx} 
\graphicspath{{S-W}} 

\newcommand{\thinrule}{\noindent\rule{\linewidth}{0.1mm}}
\usepackage{caption} 
\usepackage{wrapfig} 

\usepackage{epigraph} 
\usepackage{CormorantGaramond}
\usepackage{xcolor} 
\definecolor{OMESHICHA}{RGB}{55,107,109} 
\definecolor{ONANDOCHA}{RGB}{70,93,76} 
\definecolor{MIDORI}{RGB}{34,125,81} 
\definecolor{ULB}{RGB}{0,32,96} 
\definecolor{SHIRONERI}{RGB}{252,250,242}
\definecolor{USUAO}{RGB}{145, 180, 147}
\definecolor{SHIRONEZUMI}{RGB}{189,192,186}
\definecolor{TETSU}{RGB}{38,69,61}
\definecolor{TETSUONANDO}{RGB}{37,83,89}
\definecolor{SHINBASHI}{RGB}{0,137,167}
\definecolor{AI}{RGB}{13,86,97}
\definecolor{CHITOSEMIDORI}{RGB}{54,86,60}
\definecolor{KleinBlue}{RGB}{30, 50, 190} 
\usepackage{cite} 
\usepackage[title]{appendix} 
\usepackage{url} 
\usepackage[colorlinks, citecolor=KleinBlue, linkcolor=KleinBlue]{hyperref} 
\hypersetup{linkcolor=KleinBlue, urlcolor=KleinBlue} 
\usepackage[most]{tcolorbox} 
\usepackage{mdframed} 
\usepackage[shortlabels]{enumitem}
\usepackage{diagbox} 
\usepackage{hhline} 
\usepackage{backref}

\usepackage{tocloft}
\begin{document}
\newtheorem{theorem}{Theorem}[section]
\newtheorem{proposition}{Proposition}[section]
\newtheorem{example}{Example}[section]
\newtheorem{definition}{Definition}[section]
\newtheorem{corollary}{Corollary}[section]
\newtheorem{lemma}{Lemma}[section]
\newtheorem{remark}{Remark}[section]
\maketitle

\noindent\thinrule

\begingroup
\renewcommand{\thefootnote}{\arabic{footnote}}
\footnotetext[1]{\emph{Affiliation where this work was primarily written}: School of Mathematical Science, Capital Normal University.}
\footnotetext[2]{\emph{Current Affiliation}: Department of Mathematics, Université Libre de Bruxelles.}
\endgroup

\begingroup
\renewcommand{\thefootnote}{\emailsymbol}
\footnotetext{e-mail: \texttt{zhiyuan.liu@ulb.be}}
\endgroup

\begin{flushleft}
\textbf{\textsc{Abstract}}\par
\vspace{0.4em}

We study the topology of the regular loci of two complexified Hamiltonian integrable systems using the Zariski-van Kampen method. In particular, we show that the fundamental group of the regular locus for the  complexified planar Kepler problem is $\mathbf{F}_2\rtimes\mathbb{Z}$ where the action of the generator of $\mathbb{Z}$ interchanges the generators of $\mathbf{F}_2$, whereas that for the complexified spherical pendulum is $A\left(\tilde{A}_3\right)$, that is the affine Artin group of type $\tilde{A}_3$.  

\vspace{1ex}


\noindent\textbf{AMS Classification Code.} 37J35; 37J38; 70H06.

\tableofcontents
\end{flushleft}

\noindent\thinrule

\section{\textsc{Introduction}}
Completely integrable Hamiltonian systems play an important role in symplectic geometry and mathematical physics. In the Liouville sense, an integrable Hamiltonian system $(M,\omega,\mu, H)$ with $n$ degrees of freedom is given by a symplectic manifold $(M,\omega)$ of dimension $2n$ together with $n$ almost everywhere independent Poisson commuting first integrals $$\mu=(\mu_1=H,...,\mu_n): M\longrightarrow B\subset \mathbb{R}^n.$$ 
 The celebrated Liouville--Arnold theorem \cite{arnold1989mathematical} asserts that, near a regular compact fiber of $\mu$, the dynamics can be described in \emph{action-angle coordinates}, and the regular fibers are Lagrangian tori. Thus, from the local point of view, an integrable system looks like a torus fibration over the set of regular values $B\setminus S$ (which will be called the \emph{regular locus} through this article).
 
 From the global point of view, however, such a fibration need not be trivial, hence the action-angle coordinates are not always defined globally. This obstruction is one of the global topological phenomena in integrable systems. A precise way to measure this obstruction is the \emph{Hamiltonian monodromy}, which was initiated by Duistermaat in \cite{duistermaat1980global}. 

To recall its definition, let $\mu:M^\circ\to B\setminus S$ be the restriction of the integral map to the regular locus $B\setminus S$, and let $\Lambda\to B\setminus S$ be the period lattice bundle whose fiber $\Lambda_b$ over $b\in B\setminus S$ is the lattice of periods of the Liouville torus $\mu^{-1}(b)$.  By picking an Ehressmann connection, the parallel transport of the lattice along a loop $\gamma\in\pi_1(B\setminus S,b)$ will induce a well-defined automorphism of the lattice $\Lambda_b\cong \mathbb{Z}^n$, hence a representation
  $$\mathrm{Mon}: \pi_1(B\setminus S,b)\longrightarrow \mathrm{GL}_n(\mathbb{Z}).$$
This morphism $\mathrm{Mon}$ is called the Hamiltonian monodromy, its image is sometimes called the \emph{monodromy group} \cite{cushman2002sign}. The non-triviality of this map implies the nonexistence of the globally defined action-angle coordinates \cite{duistermaat1980global}.
Many examples of nontrivial Hamiltonian monodromy were found since its introduction. The first such example is the spherical pendulum, whose nontrivial monodromy was proved by Duistermaat in the same paper \cite{duistermaat1980global}. Subsequent work showed that the nontrivial Hamiltonian monodromy will exist if $\mu$ has a focus-focus singularity \cite{martynchuk2021recent}.
The significance of Hamiltonian monodromy goes well beyond its role as an obstruction to global action-angle coordinates. On the geometric side, it governs the global geometry of the torus fibration and it can be applied, for instance, to study the structure of the moduli spaces of some Higgs bundles \cite{baraglia2018monodromy, kydonakis2020monodromyrank2parabolic}.  On the physical side, Hamiltonian monodromy manifests itself as a defect in the joint spectral lattice of quantum integrable systems; this phenomenon has been observed, for example, in the quantum spherical pendulum and in molecular spectra \cite{zhilinskii2010quantum,sadovskii2006quantum}. We refer to \cite{martynchuk2021recent, zhilinskii2010quantum} for a more comprehensive survey on its recent developments.

Over the years, the notion of monodromy has been generalized in several directions \cite{guillemin1989monodromy, sadovskii2006quantum}. One natural generalization is the \emph{complex monodromy} when the real mechanical system admits a natural \emph{complexification} \cite{sun2022monodromy, cushman2001complex}. It can be realized, for instance, when the phase space $M$ is a real algebraic variety and the integral map $\mu$ is polynomial \cite{vanhaecke2001integrable}.  In that setting, the integrable system can be directly complexified into $(M_\C,\omega_\C,\mu_\C)$ where $M_\C$ is a complex symplectic manifold with the holomorphic symplectic $(2,0)$-form $\omega_\C$. The corresponding regular locus is in general the complement of an algebraic hypersurface in a complex affine space. By constructing periodic lattice on the generic fiber $\mu^{-1}(\mathbf{c})$ or its (partial) compactification, one can still talk about the complex Hamiltonian monodromy.  This change of viewpoint, however, enriches the topology of the regular locus and creates new monodromy phenomena invisible in the real case \cite{sun2022monodromy}. For this reason, understanding the topology of the regular locus is essential if one wants to determine the full monodromy group, rather than merely compute monodromy along a single loop.

The present paper is primarily concerned with the topology of the regular loci arising from certain complexified Hamiltonian systems. More precisely, we compute the fundamental groups of the regular loci arising from two such systems, namely the complexified spherical pendulum \cite{cushman2001complex} and the complexified planar Kepler problem\cite{sun2022monodromy}.  In the classical setting, the discriminant set $S$ merely consists of isolated points, that makes the fundamental group of the regular locus easy to determine. By contrast, in the complexified setting, the discriminant set is given by an affine algebraic curve $\mathscr{C}\subset \mathbb{C}^2$, and the corresponding regular locus turns out to be the complement $\mathbb{C}^2\setminus \mathscr{C}$, whose fundamental group is substantially more difficult to determine. Fortunately, there is an effective method in algebraic geometry, developed by Zariski and van Kampen \cite{van1933fundamental, zariski1929problem}, to attack this problem. For the complexified spherical pendulum, the existence of nontrivial monodromy was analyzed by Beukers and Cushman \cite{cushman2001complex}.  More recently, Sun and You showed that the complexified planar Kepler problem also has nontrivial monodromy and formulated a conjectural picture for the topology of the regular locus as well as the monodromy group \cite[\S6]{sun2022monodromy}. Our computations therefore also provide a description of the monodromy groups associated to this system.

The rest of this paper is organized as follows.  Section \ref{sec2} is devoted to a brief but self-contained introduction to the Zariski--van Kampen method.   We include this discussion because, in the complexified two-degree-of-freedom examples considered in this paper, the problem of determining the topology of the regular locus naturally reduces to the topology of the complement $\mathbb{C}^2\setminus\mathscr{C}$  of an affine plane curve. For this reason, the Zariski--van Kampen method should be regarded  as a general method in this class of problems. In Section \ref{sec 3}, we turn to the complexified planar Kepler problem. After reviewing its  complexification, we introduce the method in the determination of the regular locus. We then apply the machinery introduced in Section~\ref{sec2} to compute the fundamental group of the regular locus $\pi_1(\mathbb{C}^2\setminus \mathscr{C})\cong \mathbf{F}_2\rtimes\mathbb{Z}$ (Theorem \ref{theorem 3.1}).  Section~\ref{sec4} is devoted to the complexified spherical pendulum. In this case the discriminant set is more complicated: it is given by an irreducible quintic curve, and the fundamental group turns out to be more compicated. The main result of this section is that $\pi_1(\mathbb{C}^2\setminus \mathscr{C})\cong A\left(\tilde{A}_3\right)$ (Theorem \ref{theorem 4.1}), which is known as the \emph{affine Artin group} \cite{van1983homotopy} of type $\tilde{A}_3$.

\noindent\textbf{\textsc{Acknowledgement.}} The author is grateful to his master's supervisor, Prof. Shanzhong Sun, for valuable discussions and to Prof. José I. Cogolludo-Agustín for his kind and patient explanations of the Zariski--van Kampen method. Also thanks to Hongjie Zhou and Yunpeng Meng for the helpful discussion.

\section{\textsc{Braid Monodromy and Zariski-van Kampen Method}}\label{sec2}

In this section we briefly review the Zariski-van Kampen method. We refer to \cite{cogolludo2011braid} for a more detailed introduction. 

\subsection{Braid Groups}

\begin{figure}[htbp] 
    \centering
\begin{tikzpicture}[scale=0.8, thick]
    \def\px{2.5}  
    \def\py{1.5}  
    \def\pw{8}    
    \def\h{4}     

    \filldraw[fill=cyan!10, draw=cyan!40, opacity=0.8] 
        (0,0) -- (\pw,0) -- (\pw+\px,\py) -- (\px,\py) -- cycle;
    \node[text=black] at (0.5+\px/4, \py/2) {$\Pi_1$};

    \begin{scope}[shift={(0,\h)}]
        \filldraw[fill=cyan!10, draw=cyan!40, opacity=0.8] 
            (0,0) -- (\pw,0) -- (\pw+\px,\py) -- (\px,\py) -- cycle;
        \node[text=black] at (0.5+\px/4, \py/2) {$\Pi_2$};
    \end{scope}

    \coordinate (B1) at (1.5+\px/2, \py/2);
    \coordinate (B2) at (3.0+\px/2, \py/2);
    \coordinate (Bi) at (4.5+\px/2, \py/2);
    \coordinate (Bn) at (6.5+\px/2, \py/2);

    \coordinate (T1) at (1.5+\px/2, \py/2+\h);
    \coordinate (T2) at (3.0+\px/2, \py/2+\h);
    \coordinate (Ti) at (4.5+\px/2, \py/2+\h);
    \coordinate (Tn) at (6.5+\px/2, \py/2+\h);

    \foreach \i/\label in {1/1, 2/2, i/i, n/n} {
        \fill[black] (B\i) circle (2pt) node[below=2pt] {$\label$};
        \fill[black] (T\i) circle (2pt) node[above=2pt] {$\label$};
    }
    \draw[black, thick] (T1) -- (B2);
    \draw[white, line width=4.5pt] (T2) -- (B1);
    \draw[red!70!black, thick] (T2) -- (B1);

    \draw[blue!70!black, thick] (Ti) -- (Bi);

    \draw[green!60!black, thick] (Tn) -- (Bn);

    \node at (3.75+\px/2, \py/2) {$\cdots$};
    \node at (3.75+\px/2, \py/2+\h) {$\cdots$};
    \node at (5.5+\px/2, \py/2) {$\cdots$};
    \node at (5.5+\px/2, \py/2+\h) {$\cdots$};

\end{tikzpicture}
\caption{Geometric braids.} 
    \label{fig2} 
\end{figure}

Let $\Pi_1,\Pi_2$ be two parallel planes in $\mathbb{R}^3$, in particular, we assume they both parallel to the $xOy$-plane, and $\Pi_2$ is above $\Pi_1$ in the sense that $\Pi_2$ has larger $z-$component. There are $n$ marked ordered positions on each plane $\Pi_i$, namely $1,2,...,n$. A \emph{geometric braid} \cite{hansen1989braids} on $n$ strings is a family of pairwise disjoint simple arcs $(\beta_1,\beta_2,...,\beta_n)\subset\R^3$ such that each arc $\beta$ joins the $i-$th position in $\Pi_2$ to the $\sigma(i)-$th position in $\Pi_1$ for some permutation $\sigma\in S_n$, and each two arcs do not intersect,  see figure \ref{fig2}.

 If every string connects corresponding points with the same label directly, the braid is then called the \emph{identity}, denoted by $1$. Two $n-$string geometric braids $(\beta_1,...,\beta_n)$ and $(\gamma_1,...,\gamma_n)$ are called equivalent, if they can be continuously deformed to each other. The \emph{product} or \emph{composition} of two braids is defined to be the juxtaposition. Clearly, the braids product satisfy the associative law, and any braid product with the identity is itself. Hence the collection of all braids with $n$ strings forms a group, called \emph{braid groups}, denoted by $B_n$.

Although braids can be complicated, they can be written as the products of a sequel of simple braids, denoted $\sigma_i$. $\sigma_i$ is the braid with just the $i-$th and the $(i+1)-$th position interchange and only once, see figure \ref{fig4}. These $\sigma_i$ forms the generators of the braid group $B_n$, and the presentation is given by \cite{hansen1989braids}:
$$B_n=\left\langle\sigma_1,...,\sigma_{n-1}\left|\begin{cases}\sigma_i\sigma_j=\sigma_j\sigma_i,&|i-j|\geqslant 2\\\sigma_i\sigma_{i+1}\sigma_i=\sigma_{i+1}\sigma_i\sigma_{i+1},&1\leqslant i\leqslant n\end{cases}\right\rangle\right..\\$$

\begin{figure}[htbp]
    \centering
\begin{tikzpicture}[
    scale=1.0,
    line cap=round,
    line join=round,
    font=\small
]
    \colorlet{s1}{blue!55!black}
    \colorlet{s2}{teal!55!black}
    \colorlet{neutral}{black!55}
    \colorlet{ptc}{black!55}

    \draw[neutral, line width=1.05pt] (1,3) -- (1,0);
    \draw[neutral, line width=1.05pt] (4,3) -- (4,0);

    \draw[s1, line width=1.1pt]
        (2,3) .. controls (2,1.52) and (3,1.48) .. (3,0);

    \draw[white, line width=2.8pt]
        (3,3) .. controls (3,1.52) and (2,1.48) .. (2,0);
    \draw[s2, line width=1.1pt]
        (3,3) .. controls (3,1.52) and (2,1.48) .. (2,0);

    \foreach \i in {1,2,3,4} {
        \fill[ptc] (\i,3) circle (1.7pt);
        \fill[ptc] (\i,0) circle (1.7pt);
    }

    \node[above=2.5pt] at (2,3) {$i$};
    \node[above=2.5pt] at (3,3) {$i+1$};
    \node[below=2.5pt] at (2,0) {$i$};
    \node[below=2.5pt] at (3,0) {$i+1$};
\end{tikzpicture}
\caption{Simple braids $\sigma_i$} 
    \label{fig4}
\end{figure}

We now recall several standard realizations of braid groups that will be used later.

\begin{example}\label{example 2.1}
\emph{Suppose there are $n$ unordered points moving on the complex plane $\mathbb{C}$ without collisions, the motion of these $n$ points forms a configuration space $X$, that is 
$$\begin{aligned}X&=\frac{\left\{(z_1,...,z_n)\in\mathbb{C}^n\bigg|z_i\neq z_j, i\neq j\right\}}{S_n}\\
&:=X'/S_n.
\end{aligned}$$ 
The fundamental group of this configuration space is actually the braid group $$\pi_1(X,x_0)\cong B_n.$$
In fact, note that the quotient $X'\longrightarrow X$ defines an $S_n$-bundle over $X$, if we choose a section $s$ such that $s(x_0)=\tilde{x}_0\in X'$, then by the homotopy lifting property, each loop in $\pi_1(X,x_0)$ can be lifted to a path in $X'$ starting at $\tilde{x}_0$. Such a path describes a motion of $n$ ordered points in the complex plane without collisions, which is precisely an $n$-string braid. Hence $\pi_1(X,x_0)$ is actually the braid group $B_n$. $\clubsuit$}
\end{example}

\begin{example}[Artin representation]\label{example 2.2}
\emph{Another interesting model is the mapping class group $\mathscr{M}_{0,1}^n$ of the compact Riemann surface with 0 genus, 1 boundary and $n$ marked points, that is the $n-$punctured disk $\mathbb{D}\setminus\{a_1,...,a_n\}$, it was proved in \cite{birman1969mapping} that 
$$\mathscr{M}_{0,1}^n\cong B_n.$$
Mapping class group has a natural action on the fundamental group 
$$\pi_1(\mathbb{D}\setminus\{a_1,...,a_n\},a_0)=\underbrace{\mathbb{Z}*\cdots*\mathbb{Z}}_{n-\text{terms}}:=\mathbf{F}_n.$$
If we order the generators of $\mathbf{F}_n$ in an appropriate way, say $g_1,...,g_n$, then the action coincides with the Artin representation \cite{birman1969mapping}:
\begin{equation}\label{equation 1}
\begin{aligned}
\alpha: \mathscr{M}^n_{0,1}\cong B_n&\longrightarrow \mathrm{Aut}\left(\mathbf{F}_{n}\right)\\
\alpha(\sigma_i)(g_j)&=\begin{cases}g_{i+1}&j=i\\g_jg_ig_j^{-1}& j=i+1\\ g_j&\text{otherwise}\end{cases}.
\end{aligned}
\end{equation}
This example will be used in the construction of the braid monodromy. $\clubsuit$}
\end{example}

\begin{example}\label{example 2.3}
\emph{Let $\mathbb{C}[y]_n$ be the set of monic polynomials with degree $n$, let $\Delta\subset \mathbb{C}[y]_n$ be the subset consisting of polynomials with multiple roots, i.e those monic polynomials $f$ with vanishing discriminants $\Delta_f=0$, then the fundamental group $\pi_1(\mathbb{C}[y]_n\setminus\Delta)$ is exactly $B_n$.}

\emph{Indeed, by taking coefficients, the set $\mathbb{C}[y]_n\setminus\Delta$ can be identified with the set of $n-$tuples in $\mathbb{C}^n$ with non-zero discriminants, and this is equivalent to the set of $n$ distinct complex roots, it then becomes the model stated in example \ref{example 2.1}, hence its fundamental group is obviously $B_n$. $\clubsuit$}
\end{example}

\subsection{The Fundamental Group of the Complement of an Algebraic Curve}

Let $$f(x,y)=y^n+\sum_{i=1}^na_i(x)y^{n-i}\in\mathbb{C}[y]_n$$
be a monic degree $n$ polynomial, let $\mathscr{C}$ be the algebraic curve 
$$\left\{(x,y)\in\mathbb{C}^2\bigg|f(x,y)=0\right\}\subset\mathbb{C}^2.$$
We will introduce Zariski-van Kampen method to compute the fundamental group $\pi_1(\mathbb{C}^2\setminus\mathscr{C})$.

We first define the projection onto the first component 
$$\begin{aligned}p:\mathbb{C}^2\setminus \mathscr{C}&\longrightarrow \mathbb{C}\\
(x,y)&\mapsto x.\end{aligned}$$
For a generic $x$, the fiber of the projection $p$ is $\mathbb{C}\setminus\{n\,\text{points}\}$. The exceptional values are those for which the polynomial $f(x,y)=f_x(y)$ has multiple roots. These points are precisely the branch points of the projection or the singularities of the curve $\mathscr{C}$. Let $S=\{x_1,...,x_s\}\subset\mathbb{C}$ be the set of all such points. Denoted by $L_k=p^{-1}(x_k)$, and let $$\mathscr{L}=\bigcup_{k=1}^sL_k$$
be the union of all such lines, they will be called the \emph{singular fibers}, see figure \ref{figure5}.

\begin{figure}[htbp]
    \centering
\begin{tikzpicture}[line cap=round,line join=round,scale=0.8]

\coordinate (A) at (0,0);
\coordinate (B) at (8.2,0);
\coordinate (C) at (11.4,2.0);
\coordinate (D) at (3.2,2.0);
\draw[line width=1.0pt] (A)--(B)--(C)--(D)--cycle;

\fill (4.0,1.25) circle (2.2pt);
\fill (6.8,1.18) circle (2.2pt);

\node[below] at (4.0,1.05) {$x_1$};
\node[below] at (6.8,0.98) {$x_2$};
\node at (1.3,0.65) {$\mathbb{C}$};

\draw[line width=1.0pt,->] (1.3,5.2) -- (1.3,3.0);
\node[left] at (1.15,4.1) {$p$};

\draw[red,line width=0.9pt] (4.0,1.25) -- (4.0,7.55);
\draw[red,line width=0.9pt] (6.8,1.18) -- (6.8,8.25);

\node at (4.05,8.0) {$L_1$};
\node at (7.15,8.0) {$L_2$};

\fill (4.0,5.35) circle (2.2pt);
\fill (6.8,5.95) circle (2.2pt);

\draw[line width=1.0pt]
(4.0,5.35)
.. controls (4.15,6.25) and (5.25,6.45) .. (6.8,5.95)
.. controls (8.0,6.7) and (9.4,7.7) .. (10.0,8.15);

\draw[line width=1.0pt]
(4.0,5.35)
.. controls (4.25,4.55) and (5.0,4.35) .. (6.8,5.95)
.. controls (8.2,5.55) and (9.4,5.2) .. (10.1,5.05);

\node at (8.75,7.55) {$\mathscr{C}$};

\end{tikzpicture}
\caption{projection onto the first component}
\label{figure5}
\end{figure}

After a  suitable change of coordinates, we may assume that the projection $p$ is \emph{generic} with respect to $\mathscr C$ in the sense that the curve $\mathscr C$ has no vertical component, distinct critical points of the projection have distinct $x$-coordinates, and for each critical value $x_k$, the fiber $L_k=p^{-1}(x_k)$ contains exactly one non-transversal point of $L_k\cap \mathscr C$, while all other intersections are transverse. Since an affine automorphism of $\C^2$ induces a homeomorphism of $\C^2\setminus \mathscr C$, this assumption does not affect the computation of $\pi_1(\C^2\setminus \mathscr C)$.

\begin{proposition}[\cite{cohen1997braid}]
The restriction of the projection $p$:
$$p:\mathbb{C}^2\setminus(\mathscr{C}\cup\mathscr{L})\longrightarrow \mathbb{C}\setminus S$$
is a fibration with fiber $F=\mathbb{C}\setminus \{n \,\text{points}\}$. The structure group of this fiber bundle is precisely the braid group $B_n$.
\end{proposition}

Choose a base point $x_0\in \mathbb{C}\setminus S$ and $y_0\in F$, the fundamental groups of the base manifold and the fiber are simply:
$$\pi_1(\mathbb{C}\setminus S,x_0)\cong \mathbf{F}_{s},\,\pi_1\left(F,y_0\right)\cong\mathbf{F}_{n}.$$ 
The fundamental group of the base manifold has an action on the fundamental group of the fiber, this action is called \textbf{braid monodromy}, denoted by $\rho$:
$$\rho: \mathbf{F}_{s}\longrightarrow\mathrm{Aut}\left(\mathbf{F}_{n}\right).$$
The construction of the braid monodromy $\rho$ is as follows.

We take the algebraic curve $f(x,y)=0$ as a map:
$$\begin{aligned}f: \mathbb{C}&\longrightarrow\mathbb{C}[y]_n\cong\mathbb{C}^n\\
 x&\mapsto f_x(y)=f(x,y).\end{aligned}$$
Let $\Delta$ be the subset containing the polynomials with vanishing discriminants (see example \ref{example 2.3}), observe that $f^{-1}(\Delta)$ is precisely $S$, hence the restriction of $f$ on $\mathbb{C}\setminus S$ induces a homomorphism of fundamental groups:
$$\begin{tikzcd}
{f_*: \pi_1(\mathbb{C}\setminus S,x_0)} \arrow[d, Rightarrow, no head, shift left] \arrow[r] & \pi_1(\mathbb{C}^n\setminus\Delta) \arrow[d, Rightarrow, no head] \\
\mathbf{F}_s                                                                                    & B_n                                                              
\end{tikzcd}$$
The braid monodromy is given by the composition of Artin representation $\alpha$ (see example \ref{example 2.2}) and $f_*$:
$$\rho=\alpha\circ f_*:  \mathbf{F}_{s}\longrightarrow\mathrm{Aut}\left(\mathbf{F}_{n}\right).$$ 
We can now apply homotopy exact sequence of the fibration:
$$\begin{tikzcd}
1 \arrow[r] & {\pi_1(F,y_0)} \arrow[r] \arrow[d, Rightarrow, no head] & \pi_1(\mathbb{C}^2\setminus\mathscr{C}) \arrow[r, "p_*"] & {\pi_1(\mathbb{C}\setminus S,x_0)} \arrow[r] \arrow[d, Rightarrow, no head] & 1. \\
            & \mathbf{F}_n                                                                &                                                          & \mathbf{F}_s                                            &  
\end{tikzcd}$$
Denoted by $\gamma_1,...,\gamma_s$, $g_1,...,g_s$ the generators of the fundamental groups $\pi_1(\mathbb{C}\setminus S,x_0)$ and $\pi_1(F,y_0)$ respectively (see figure \ref{figure6}), then the main result of the Zariski-van Kampen can be stated as follows:

\begin{figure}[htbp]
\centering
\includegraphics[width=0.7\linewidth]{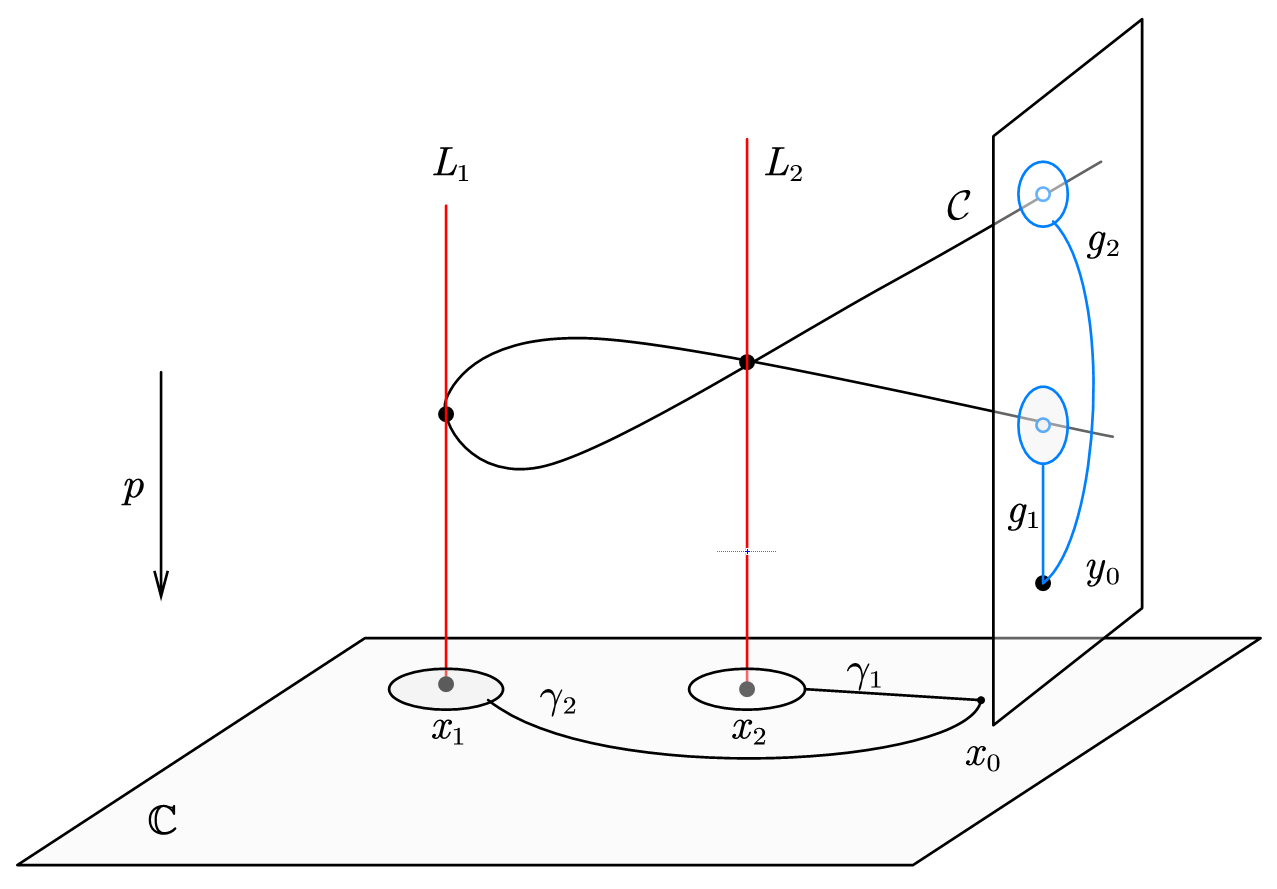}
\caption{generators}\label{figure6}
\end{figure}

\begin{theorem}[Zariski-van Kampen, \cite{cogolludo2011braid}]\label{theorem 2.2}
The fundamental group of total space is the semi-product along braid monodromy:
$$\pi_1\left(\mathbb{C}^2\setminus(\mathscr{C}\cup\mathscr{L})\right)\cong \mathbf{F}_{s}\ltimes_{\rho}\mathbf{F}_{n}$$
In particular, it has the presentation
$$\pi_1\left(\mathbb{C}^2\setminus(\mathscr{C}\cup\mathscr{L})\right)=\left\langle\gamma_1,...,\gamma_s,g_1,...,g_n\bigg|\gamma_k^{-1}g_i\gamma_k=\rho(\gamma_k)(g_i)\right\rangle.$$
Moreover, the fundamental group of $\mathbb{C}^2\setminus\mathscr{C}$ is the quotient by the normal closure generated by $\gamma_1,...,\gamma_s$, and it has the presentation
$$\pi_1\big(\mathbb{C}^2\setminus\mathscr{C}\big)=\left\langle g_1,...,g_n\bigg|g_i=\rho(\gamma_k)(g_i) \right\rangle.$$
\end{theorem}

\begin{example}[The Riemann surface of $\sqrt{z}$]\label{example 2.4}
\emph{As a simple example, let's consider the fundamental group of the complement of the curve $X_{\sqrt{z}}:y^2=x$, which is the Riemann surface of the square root function.}

\emph{Observe that the curve $X_{\sqrt{z}}$ has only one branched point $(0,0)$ in projecting onto the first component, hence $\pi_1(\mathbb{C}\setminus S)\cong \mathbb{Z}$, and braid monodromy will be given in $B_2$. Now choose $\gamma(t)=e^{2\pi\sqrt{-1}}t$ be the generator in $\pi_1(\mathbb{C}\setminus S)$, observe that the change of $y$ is simply $\pm e^{\pi\sqrt{-1}t}$, hence it yields the braid as $\sigma_1$, by (\ref{equation 1}), the generators $g_1,g_2$ on the fiber should satisfy the relation $g_1=g_2$, hence $$\pi_1\left(\mathbb{C}^2\setminus X_{\sqrt{z}}\right)\cong\mathbb{Z}.$$}
$\clubsuit$
\end{example}

We refer to \cite{cogolludo2011braid} for more interesting examples. 

\section{\textsc{The Complexified Planar Kepler Problem}}\label{sec 3}

Let's first recall the classical planar Kepler problem, for more detailed contents can be found in \cite{cordani2003separation}.

The phase space of the planar Kepler problem is the symplectic manifold $T^*(\mathbb{R}^2\setminus\{\mathbf{0}\})\cong(\mathbb{R}^2\setminus\{\mathbf{0}\})\times\mathbb{R}^2$ endowed with the standard symplectic form:
$$\omega=dp_1\wedge dq_1+dp_2\wedge dq_2$$
where $(q_1,q_2,p_1,p_2):=(\mathbf{q,p})$ denotes the coordinates on $T^*(\mathbb{R}^2\setminus\{\mathbf{0}\})$. The Hamiltonian $H$ is defined by
\begin{equation}\label{equation 5}H(\mathbf{q,p})=\frac{1}{2}\left(p_1^2+p_2^2\right)-\frac{1}{\sqrt{q_1^2+q_2^2}}.\end{equation}

This system carries a natural Hamiltonian action of $\mathrm{SO}(2)$, the associated moment map is the angular moment
$$\begin{aligned}J: T^*\left(\mathbb{R}^2\setminus\{\mathbf{0}\}\right)&\longrightarrow\mathbb{R}\\ J(\mathbf{q,p})&=q_1p_2-q_2p_1\end{aligned}$$
hence the system admits two Poisson commuting independent first integrals, namely $H$ and $J$, and it is therefore Liouville integrable. 

\subsection{The Complexification}

We cannot simply replace the configuration space $\mathbb{R}^2\setminus\{\mathbf{0}\}$ by $\mathbb{C}^2\setminus\{\mathbf{0}\}$ in the complexification, since the square root term in (\ref{equation 5}) becomes multi-valued. To avoid this ambiguity, we define the complexified configuration space to be the punctured algebraic surface
$$Q=\left\{ (w,z_1,z_2)\in\mathbb{C}^3\bigg|w^2=z_1^2+z_2^2\right\}\setminus\{\mathbf{0}\}.$$

\begin{proposition}
The complexified configuration space $Q$ is a $2$ dimensional complex symplectic manifold with trivial cotangent bundle.
\end{proposition}

\noindent\emph{Proof}. Note that the algebraic surface $w^2=z_1^2+z_2^2$ has the only singularity at the origin, hence $Q$ is a holomorphic complex 2-dimensional manifold. 

There are 2 holomorphic charts $(Q_1,\varphi_1)$ and $(Q_2,\varphi_2)$ on $Q$, namely 
\begin{equation}\label{equation 6}
\left\{\begin{aligned}
&Q_1=\left\{(w,z_1,z_2)\in Q\bigg| w\in\mathbb{C}\setminus\mathbb{R}_{>0}\right\}\\
&\varphi_1(w,z_1,z_2)=(z_1,z_2)\\
\end{aligned}\right.\end{equation}
and
\begin{equation}
\left\{\begin{aligned}&Q_2=\left\{(w,z_1,z_2)\in Q\bigg| w\in\mathbb{C}\setminus\mathbb{R}_{<0}\right\}\\
&\varphi_2(w,z_1,z_2)=(z_1,z_2).
\end{aligned}\right.
\end{equation}
If we write $z_1^2+z_2^2=re^{\sqrt{-1}\theta}$, where $\theta\in [-\pi,\pi)$, then the inverses on each coordinate chart are respectively given by
$$\begin{aligned}\varphi_1^{-1}(z_1,z_2)&=\left(\sqrt{r}e^{\frac{\sqrt{-1}}{2}(\theta+2\pi)},z_1,z_2\right)\\
\varphi_2^{-1}(z_1,z_2)&=\left(\sqrt{r}e^{\frac{\sqrt{-1}}{2}\theta},z_1,z_2\right).\end{aligned}$$

Clearly, the transition $\varphi_2\circ\varphi_1^{-1}$ on the overlap is the identity, thus $Q$ has trivial tangent and cotangent bundle. Consequently, the complexified phase space is simply $T^*Q=Q\times\mathbb{C}^2$, the local coordinate will be denoted by $(z_1,z_2,w_1,w_2):=(\mathbf{z,w})$. Under these notations, the standard symplectic form on $T^*Q$ can be written as
$$\omega=dz_1\wedge dw_1+dz_2\wedge dw_2.$$
In particular, $\omega$ is a holomorphic $(2,0)-$form on $T^*Q$, hence a holomorphic symplectic manifold. $\clubsuit$

Let $(z_1,z_2)$ be the coordinate on $Q_i$, the complexified Hamiltonian is given by 
$$H(\mathbf{z,w})=\frac{1}{2}\left(w_1^2+w_2^2\right)-\frac{(-1)^i}{\sqrt{z_1^2+z_2^2}},\,\,i=1,2.$$
It now becomes a holomorphic single-valued function in $\mathscr{O}(T^*Q)$. The complex Lie group $\mathrm{SO}(2,\mathbb{C})$ acts Hamiltonianly on $T^*Q$ by

$$\begin{pmatrix}a&-b\\b&a\end{pmatrix}\cdot(z_1,z_2,w_1,w_2):=\left((z_1,z_2)\begin{pmatrix}a&-b\\b&a\end{pmatrix}, (w_1,w_2)\begin{pmatrix}a&-b\\b&a\end{pmatrix}\right)$$ 
and the moment map associated to this action is the complexified angular moment:
$$\begin{aligned}J: T^*Q&\longrightarrow \mathbb{C}\\
J(\mathbf{z,w})&=z_1w_2-z_2w_1.\end{aligned}$$
Similar to the real case, the system $(T^*Q,\omega, H, J)$ becomes a complex integrable system. The holomorphic integral map is
$$\begin{aligned}\mu: T^*Q&\longrightarrow \mathbb{C}^2\\ (\mathbf{z,w})&\mapsto (H(\mathbf{z,w}), J(\mathbf{z,w})).\end{aligned}$$ 

\subsection{Determine the Regular Locus}
We shall first determine the discriminant set of the integral map $\mu: T^*Q\longrightarrow\mathbb{C}^2$. Recall that in the complexified setting, a value $\mathbf{c}\in\mathbb{C}^2$ is said to lie in the discriminant set, if the fiber $\mu^{-1}(\mathbf{c})$ is either non-smooth or on each open neighborhood $U_\mathbf{c}$ of $\mathbf{c}$, there is a dense subset $V\subset U_\mathbf{c}$ such that the topology of $\mu^{-1}(v)$ is different from $\mu^{-1}(\mathbf{c})$ for all $v\in V$. Although the result was established in \cite[Proposition 3.1]{sun2022monodromy}, it is still worthwhile to briefly recall it here as it illustrates a standard method for determining the discriminant set in this type of problem.

On each coordinate chart $(Q_i,\varphi_i)$ of $Q$ (defined in (\ref{equation 6})), we define the following variables: 
\begin{equation}\label{eq2}
\begin{aligned}
w&:=(-1)^i\sqrt{z_1^2+z_2^2},\,\,z:=z_1w_1+z_2w_2+\sqrt{-1}(z_2w_1-z_1w_2)\\
x&:=z_1w_2-z_2w_1,\,\,\,\,\,\,\,\,\, y:=\frac{1}{2}\left(w_1^2+w_2^2\right)-\frac{1}{w}.
\end{aligned}
\end{equation}
Note that the new variables $x,y$ are just angular moment and total energy respectively. The fiber of the integral map $\mu$ has the following description:

\begin{proposition}[\cite{sun2022monodromy}]\label{proposition 4.2} Under the notation above, the discriminant set of the complexified planar Kepler problem is the algebraic curve
$$\mathscr{C}=\left\{(x,y)\in\mathbb{C}^2\bigg|y\big(1+2x^2y\big)=0\right\}.$$
\end{proposition} 

\noindent\emph{Sketch of the proof}. The key idea is to notice that the $\mathrm{SO}(2,\C)\cong\mathbb{C}^*$-action on the fiber $\mu^{-1}(\mathbf{c})$ is free and proper \cite[Prop.3.1.(1)]{sun2022monodromy}, thus $\mu^{-1}(\mathbf{c})$ can be viewed as a $\mathbb{C}^*-$bundle over the quotient $\mathscr{N}_{\mathbf{c}}:=\mu^{-1}(\mathbf{c})/\mathrm{SO}(2,\C)$. For each value $\mathbf{c}=(x,y)\in\mathbb{C}^2$, the quotient  $\mathscr{N}_{\mathbf{c}}$ can be described by 
\begin{equation}\label{equation 8}
\left\{(w,z)\in\mathbb{C}^2\bigg|2yw^2+2w=z^2+2\sqrt{-1}xz\right\}\setminus\left\{(0,0),\left(0,-2\sqrt{-1}x\right)\right\}.
\end{equation}
where $w,z$ are defined in (\ref{eq2}) \cite[Prop.3.1.(3)]{sun2022monodromy}. Therefore, to determine the discriminant set, it suffices to analyze the topology of $\mathscr{N}_\mathbf{c}$. 

Let $\mathbf{c}=(x,y)$ be fixed, let $\mathscr{R}$ be the curve defined by the polynomial 

\begin{equation}\label{eqR}p(z,w)=2yw^2+2w-z^2-2\sqrt{-1}xz\in\mathbb{C}[z,w].\end{equation}

We can now discuss in cases:
\begin{enumerate}[(i), itemsep=-1pt]
    \item  If $y=0$, then $\mathscr{R}$ is 
$$2w-z^2-2\sqrt{-1}xz=2w-x^2-\left(z+\sqrt{-1}x\right)^2=0,$$ 
the curve is isomorphic to $\mathbb{C}$.

\item If $y\neq 0$, but $\Delta=2(1+2x^2y)=0$, then $\mathscr{R}$ becomes to 
$$2y\left(w+\frac{1}{2y}\right)^2-\left(z+\sqrt{-1}x\right)^2=0,$$
the curve is isomorphic to a singular cone, and this case correspond to the singular fiber.

\item If $y\neq 0$ and $\Delta=2(1+2x^2y)\neq0$, then $\mathscr{R}$ becomes to 
$$2y(w-r)(w-s)-\left(z+\sqrt{-1}x\right)^2=0$$
for some $r\neq s$, then the curve $\mathscr{R}$ is isomorphic to $\mathbb{C}^*$. This is the case such that the fiber $M_{\mathbf{c}}$ is a generic fiber. 
\end{enumerate}

\noindent To sum up, the discriminant set of $\mu$ is indeed the curve $y(1+2x^2y)=0$. $\clubsuit$

It is straightforward to see that the integral map $\mu$ is surjective onto $\mathbb{C}^2$. Indeed, for any $(x,y)\in\mathbb{C}^2$, the defining equation (\ref{equation 8}) is a quadratic in $  w  $ and always admits solutions in $\mathbb{C}$, so the fiber $\mu^{-1}(x,y)$ is non-empty. Hence the regular locus is precisely $\mathbb{C}^2\setminus\mathscr{C}$.

\subsection{The Topology of the Regular Locus}

Now we can state and prove the main result of this section:

\begin{theorem}\label{theorem 3.1}
The fundamental group of the regular locus  $\mathbb{C}^2\setminus\mathscr{C}$ is
$$\pi_1\left(\mathbb{C}^2\setminus\mathscr{C}\right)\cong \mathbf{F}_2\rtimes\mathbb{Z},$$
where the generator of $\mathbb{Z}$ interchanges the generators of $\mathbf{F}_2$.
\end{theorem}

\noindent\emph{Proof}. Let $\mathscr{C}_1$ be the curve $1+2x^2y=0$ and $\ell$ the asymptotic line $y=0$. We regard $\mathbb{C}^2\setminus(\mathscr{C}_1\cup\ell)$ as a subset in the projective plane $\mathbb{P}^2$ via
$$(x,y)\mapsto [x:y:1]\in\mathbb{P}^2.$$
Then the line $z=0$ becomes the line at infinity of $\mathbb{C}^2$, which we denote by $\ell_{\infty}$. We shall use $\bar{\mathscr{C}_1}$ and $\bar{\ell}$ represent for the projectivization of $\mathscr{C}_1$ and $\ell$ in $\mathbb{P}^2$ respectively. In particular, $\bar{\mathscr{C}}_1$ and $\bar{\ell}$ intersect at $[1:0:0]\in\ell_{\infty}$.

Since the projective plane $\mathbb{P}^2$ is obtained from $\mathbb{C}^2$ by adding its line at infinity, we have
$$\pi_1\left(\mathbb{C}^2\setminus(\mathscr{C}_1\cup\ell)\right)=\pi_1\left(\mathbb{P}^2\setminus\left(\bar{\mathscr{C}}_1\cup\bar{\ell}\cup\ell_{\infty}\right)\right).$$
Now, we change the coordinates in $\mathbb{P}^2$ by 
$$A=\begin{bmatrix}1&0&0\\0&0&1\\0&1&0\end{bmatrix}\in\mathrm{PGL}_3(\mathbb{C}).$$
Observe that the line $\bar{\ell}$ becomes to the line at infinity $\ell_{\infty}$ after changing coordinates by $A$. Let $\bar{\mathscr{C}}_1',\ell_{\infty}'$ be the new curve after changing the coordinates, they are 
$$\bar{\mathscr{C}}_1': y^3+2x^2z=0,\,\,\ell_{\infty}': y=0.$$
Since projective transformation doesn't change the fundamental group, we have 
$$\begin{aligned}\pi_1\left(\mathbb{P}^2\setminus\left(\bar{\mathscr{C}}_1\cup\bar{\ell}\cup\ell_{\infty}\right)\right)&\cong\pi_1\left(\mathbb{P}^2\setminus\left(\bar{\mathscr{C}}'_1\cup\ell'_{\infty}\cup\ell_{\infty}\right)\right)\\
&=\pi_1\left(\left(\mathbb{P}^2\setminus\ell_{\infty}\right)\setminus\left(\bar{\mathscr{C}}'_1\cup\ell'_{\infty}\right)\right).
\end{aligned}$$

\begin{figure}[htbp]
\centering
\includegraphics[width=0.67\linewidth]{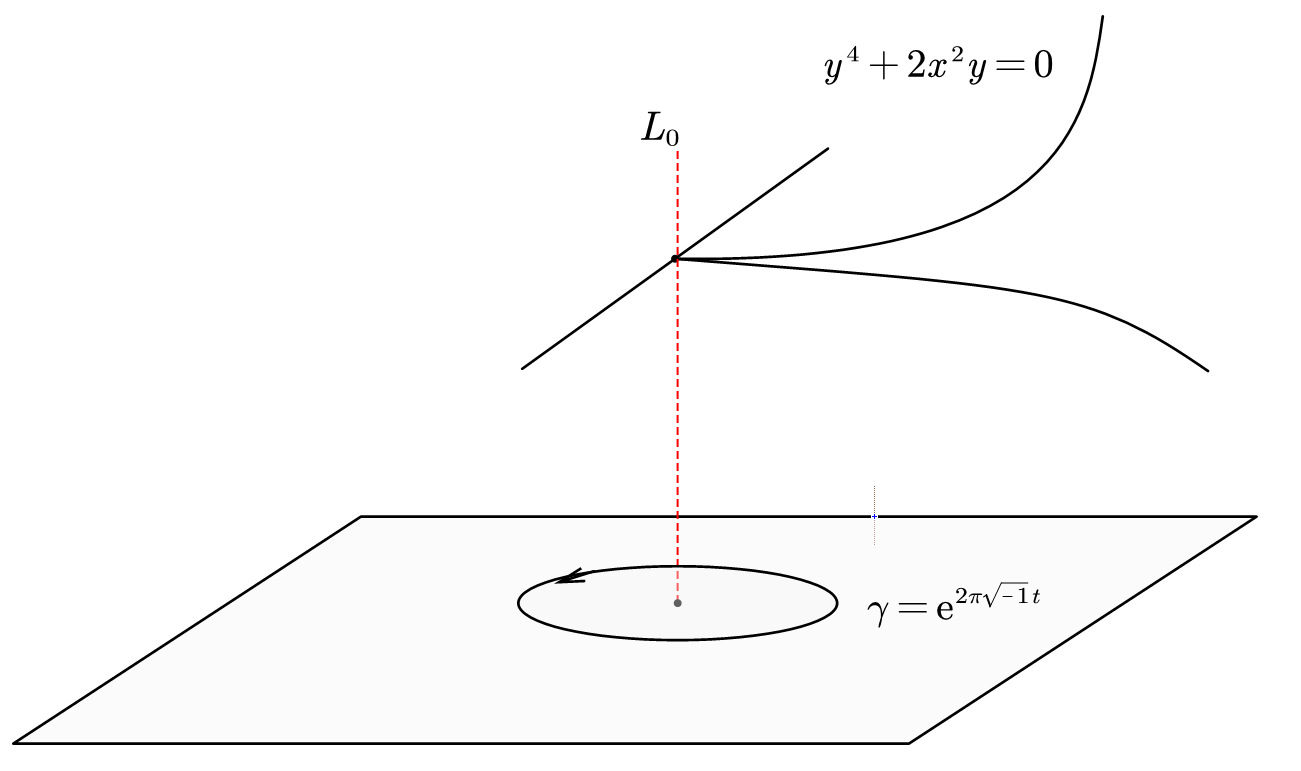}
\caption{the curve $y(2x^2+y^3)=0$}\label{figure7}
\end{figure}

Notice that $\mathbb{P}^2\setminus\ell_{\infty}$ is the affine part of $\mathbb{P}^2$, we can use the coordinate chart:
$$\begin{aligned}
\psi: \mathbb{P}^2\setminus\ell_{\infty}&\longrightarrow \mathbb{C}^2\\
[x:y:z]&\mapsto \left(\frac{x}{z},\frac{y}{z}\right)
\end{aligned}$$
so that $\left(\mathbb{P}^2\setminus\ell_{\infty}\right)\setminus\left(\bar{\mathscr{C}}'_1\cup\ell'_{\infty}\right)$ can be written under $\psi$ by
$$\left(\mathbb{P}^2\setminus\ell_{\infty}\right)\setminus\left(\bar{\mathscr{C}}'_1\cup\ell'_{\infty}\right)\cong \mathbb{C}^2\setminus\left\{y(2x^2+y^3)=0\right\}.$$

Next, we will use Zariski-van Kampen method to compute the fundamental group of the complement of the curve $y(2x^2+y^3)=0$.

The curve $y(2x^2+y^3)=0$ has only one singular point at the origin $(0,0)$, where the image of real part is depicted in figure \ref{figure7}. The singular set $S$ is therefore an isolated point and the fundamental group $\pi_1(\mathbb{C}\setminus S,x_0)$ is $\mathbb{Z}$. For each $x\in\mathbb{C}\setminus S$, since the curve has degree 4, the fundamental group of the fiber is just $$\pi_1(\mathbb{C}\setminus\{\text{4 points}\})\cong \mathbf{F}_4.$$
Their generators will be denoted by $g_1,g_2,g_3$, and $g_4$. 

Choose $\gamma(t)=e^{2\pi\sqrt{-1}t}\in\pi_1(\mathbb{C}\setminus S,x_0)$ to be the generator (see figure \ref{figure7}), then the change on the fiber will give the braid as $(\sigma_1\sigma_3\sigma_2\sigma_1)^2\in B_4$, see figure \ref{figure8}.
\begin{figure}[htbp]
    \centering
\begin{tikzpicture}[line cap=round,line join=round,>=Stealth]
\def\R{2.2}     
\def\r{0.52}    
\def\H{6.3}     
\colorlet{colone}{blue}
\colorlet{coltwo}{green}
\colorlet{colthree}{teal!70!black}
\colorlet{colfour}{red!75!black}
\draw ({-\R},0) arc[start angle=180,end angle=0,x radius=\R,y radius=\r];
\draw[dashed] ({-\R},-\H) arc[start angle=180,end angle=0,x radius=\R,y radius=\r];
\draw ({-\R},0) -- ({-\R},-\H);
\draw ({ \R},0) -- ({ \R},-\H);
\draw[line width=1.1pt,coltwo] (0,0) -- (0,-\H);

\draw[dashed,line width=1.1pt,colone,smooth,samples=60,domain=0:0.375,variable=\t]
  plot ({\R*cos(90+240*\t)},{-\H*\t + \r*sin(90+240*\t)});
\draw[line width=1.1pt,colone,smooth,samples=100,domain=0.375:1,variable=\t]
  plot ({\R*cos(90+240*\t)},{-\H*\t + \r*sin(90+240*\t)});

\draw[line width=1.1pt,colthree,smooth,samples=100,domain=0:0.625,variable=\t]
  plot ({\R*cos(210+240*\t)},{-\H*\t + \r*sin(210+240*\t)});
\draw[dashed,line width=1.1pt,colthree,smooth,samples=60,domain=0.625:1,variable=\t]
  plot ({\R*cos(210+240*\t)},{-\H*\t + \r*sin(210+240*\t)});

\draw[line width=1.1pt,colfour,smooth,samples=40,domain=0:0.125,variable=\t]
  plot ({\R*cos(330+240*\t)},{-\H*\t + \r*sin(330+240*\t)});
\draw[dashed,line width=1.1pt,colfour,smooth,samples=100,domain=0.125:0.875,variable=\t]
  plot ({\R*cos(330+240*\t)},{-\H*\t + \r*sin(330+240*\t)});
\draw[line width=1.1pt,colfour,smooth,samples=40,domain=0.875:1,variable=\t]
  plot ({\R*cos(330+240*\t)},{-\H*\t + \r*sin(330+240*\t)});

\draw ({-\R},0) arc[start angle=180,end angle=360,x radius=\R,y radius=\r];
\draw ({-\R},-\H) arc[start angle=180,end angle=360,x radius=\R,y radius=\r];
\fill[colone] ({\R*cos(90)},{\r*sin(90)}) circle (1.8pt);
\node[text=colone,above=2pt] at ({\R*cos(90)},{\r*sin(90)}) {$1$};

\fill[colthree] ({\R*cos(210)},{\r*sin(210)}) circle (1.8pt);
\node[text=colthree,left=2pt] at ({\R*cos(210)},{\r*sin(210)}) {$3$};

\fill[colfour] ({\R*cos(330)},{\r*sin(330)}) circle (1.8pt);
\node[text=colfour,right=2pt] at ({\R*cos(330)},{\r*sin(330)}) {$4$};

\fill[coltwo] (0,0) circle (1.8pt);
\node[text=coltwo,below=2pt] at (0,0) {$2$};

\fill[colthree] ({\R*cos(90)},{-\H+\r*sin(90)}) circle (1.8pt);
\node[text=colthree,above=2pt] at ({\R*cos(90)},{-\H+\r*sin(90)}) {$3$};

\fill[colfour] ({\R*cos(210)},{-\H+\r*sin(210)}) circle (1.8pt);
\node[text=colfour,left=2pt] at ({\R*cos(210)},{-\H+\r*sin(210)}) {$4$};

\fill[colone] ({\R*cos(330)},{-\H+\r*sin(330)}) circle (1.8pt);
\node[text=colone,right=2pt] at ({\R*cos(330)},{-\H+\r*sin(330)}) {$1$};
\fill[coltwo] (0,-\H) circle (1.8pt);
\node[text=coltwo,below=2pt] at (0,-\H) {$2$};

\end{tikzpicture}
\caption{The braid $\rho(\gamma)=\bigl(\sigma_1\sigma_3\sigma_2\sigma_1\bigr)^2\in B_4$.}
\label{figure8}
\end{figure}

We can now compute by (\ref{equation 1}):
\begin{equation}\label{eq relation}\begin{aligned}
g_1&=(\sigma_1\sigma_3\sigma_2\sigma_1)^2(g_1)=g_4g_3g_4^{-1},\\
g_2&=(\sigma_1\sigma_3\sigma_2\sigma_1)^2(g_2)=g_4g_3g_2g_3^{-1}g_4^{-1},\\
g_3&=(\sigma_1\sigma_3\sigma_2\sigma_1)^2(g_3)=g_4g_3g_2g_4g_2^{-1}g_3^{-1}g_4^{-1},\\
g_4&=(\sigma_1\sigma_3\sigma_2\sigma_1)^2(g_4)=g_4g_3g_2g_1g_2^{-1}g_3^{-1}g_4^{-1}.
\end{aligned}\end{equation}
Set $h=g_4g_3g_2g_1$ and $g=g_1$, By theorem \ref{theorem 2.2}, (\ref{eq relation}) therefore gives 
\begin{equation}\label{eq presentation}\begin{aligned}\pi_1\bigg(\mathbb{C}^2\setminus\left\{y(1+2x^2y)=0\right\}\bigg)&\cong\pi_1\bigg(\mathbb{C}^2\setminus\left\{y^4+2x^2y=0\right\}\bigg)\\
&=\bigg\langle g,h\bigg| h^2g=gh^2\bigg\rangle\\
&\cong\mathbf{F}_2\rtimes\mathbb{Z}.\end{aligned}\end{equation}
The theorem is proved. $\clubsuit$

As an application of Theorem \ref{theorem 3.1}, we can now determine the monodromy group of the complexified planar Kepler problem. 

First note that $g$ in (\ref{eq presentation}) can be understood as a meridian around $\mathscr{C}_1=\{1+2x^2y=0\}$, where as $h$ can be understood as the meridian around $\ell=\{y=0\}$, with suitable choices of connecting paths to the base point. For a regular value $\mathbf{c}\in\C^2\setminus\mathscr{C}$, the periodic lattice $\Lambda_{\mathbf{c}}$ of its fiber was constructed in \cite[\S4]{sun2022monodromy}, it turns out to be a rank 2 lattice in $\C^2$, and the monodromy $\mathrm{Mon}:\pi_1(\C^2\setminus\mathscr{C})\longrightarrow\mathrm{GL}_2(\Z)$ of this lattice around each generator $g,h$ was also computed in \cite[\S 5]{sun2022monodromy}, they are both conjugated to the matrix $\begin{pmatrix}-1& 0\\-2 &1\end{pmatrix}$. That is to say the normal closure of $\left\langle g^2,h^2\right\rangle$ is contained in the kernel $\ker(\mathrm{Mon})$. Therefore, by the first isomorphism theorem and the result we obtained in Theorem \ref{theorem 3.1}, we can deduce that the monodromy group is a quotient of
$$\bigg\langle g,h\bigg|gh^2=h^2g,g^2=1,h^2=1\bigg\rangle\cong\mathbb{Z}_2*\mathbb{Z}_2.$$

\section{\textsc{The Complexified Spherical Pendulum}}\label{sec4}

The configuration space of the classical spherical pendulum is the unit sphere $$S^2=\left\{(x_1,x_2,x_3)\in\mathbb{R}^3\bigg| \sum_{i=1}^3x_i^2=1\right\},$$
and its phase space is the cotangent bundle $T^*S^2$. If $\langle\cdot,\cdot\rangle$ denotes the standard Euclidean inner product on $\mathbb{R}^3$, then the phase space can be written as
$$T^*S^2=\left\{(\mathbf{x},\mathbf{v})\in S^2\times\mathbb{R}^3\bigg|\langle \mathbf{x},\mathbf{v}\rangle=0\right\}.$$
 The Hamiltonian is given by
$$H(\mathbf{x},\mathbf{v})=\frac{1}{2}\|\mathbf{v}\|^2+x_3: T^*S^2\longrightarrow\mathbb{R}.$$
As in the previous section, the Lie group $\mathrm{SO}(2)$ acts on $T^*S^2$:
$$\begin{aligned}\mathrm{SO}(2)\times T^*S^2 &\longrightarrow T^*S^2\\ (A,(\mathbf{x},\mathbf{v}))&\mapsto \left(\begin{pmatrix}A&0\\0&1\end{pmatrix}\mathbf{x}, \begin{pmatrix}A&0\\0&1\end{pmatrix}\mathbf{v}\right)\end{aligned}$$
in a Hamiltonian fashion, the moment map is the angular moment: 
$$J(\mathbf{x},\mathbf{v})=x_1v_2-x_2v_1: T^*\mathbb{S}^2\longrightarrow\mathbb{R}.$$
Together with the Hamiltonian $H$, the spherical pendulum is a Liouville integrable system.

\subsection{The Complexification}
Since everything here is polynomial, the system admits a direct complexification. Let $\langle\cdot,\cdot\rangle$ be the standard symmetric complex bilinear form on $\mathbb{C}^3$, the configuration space is now a complexified sphere:
$$S^2_\mathbb{C}=\left\{\mathbf{x}\in\mathbb{C}^3| \langle\mathbf{x},\mathbf{x}\rangle=1\right\}$$ 
and the complexified phase space is simply 
$$T^*S^2_\mathbb{C}=\left\{(\mathbf{x},\mathbf{v})\in S^2\times\mathbb{C}^3\bigg|\langle \mathbf{x},\mathbf{v}\rangle=0\right\}\subset\mathbb{C}^3\times\mathbb{C}^3$$
This is a holomorphic symplectic manifold, in particular, it is an affine variety in $\mathbb{C}^3\times\mathbb{C}^3$. 
Similarly, the complex Lie group $\mathrm{SO}(2,\mathbb{C})$ acts on $T^*S^2_\mathbb{C}$ in a Hamiltonian way and the associated moment map is
$$J(\mathbf{x},\mathbf{v})=x_1v_2-x_2v_1: T^*S^2_\mathbb{C}\longrightarrow\mathbb{C}$$ 
The complexified integral map will still be denoted by $\mu=(H,J)$. Thus the complexified spherical pendulum $(T^*S^2_\C,\omega, H,J)$ becomes a complex integrable system.

\subsection{The Regular Locus and Its Fundamental Group}

Following \cite{cushman2001complex}, we introduce the following new coordinates $$(w_1,w_2,w,z_1,z_2,z):=(\mathbf{w},\mathbf{z})$$ on $\mathbb{C}^3\times\mathbb{C}^3$:
$$\begin{aligned}w_1&=x_1+\sqrt{-1}x^2,\,\, w_2=x_1-\sqrt{-1}x_2\\z_1&=v_1+\sqrt{-1}v_2,\,\, z_2=v_1-\sqrt{-1}v_2\end{aligned}$$
In these coordinates, the defining equations for $T^*S^2_\mathbb{C}$ become
$$T^*S^2_\mathbb{C}: \begin{cases}w_1w_2+w^2=0\\ w_1z_2+w_2z_1+2wz=0.\end{cases}$$
The integral map can then be written as:
$$\mu(\mathbf{w},\mathbf{z})=(x,y):=\left(\frac{z_1z_2+z^2}{2}+w,\frac{w_2z_1-w_1z_2}{\sqrt{-1}}\right).$$ 

\begin{proposition}[\cite{cushman2001complex}]\label{proposition 3.1}
Under the notation above, the discriminant set of the complexified spherical pendulum is
\begin{equation}\label{equation 3}\mathscr{C}=\left\{(x,y)\in\mathbb{C}^2\left|\frac{27}{4}y^4+2xy^2(x^2-9)-4(x^2-1)^2=0\right.\right\}.\end{equation}
\end{proposition}
The proof follows the same strategy as in Proposition \ref{proposition 4.2} and can be found in detail in \cite[Proposition 2.1]{cushman2001complex}. It is also straightforward to check that $\mu$ is surjective and hence the regular locus is precisely $\C^2\setminus \mathscr{C}$.

\begin{theorem}\label{theorem 4.1}
The fundamental group of the regular locus of the complexified spherical pendulum is $$
\pi_1\bigg(\mathbb C^2\setminus\mathscr C\bigg)
\cong
\left\langle
a_1,a_2,a_3,a_4
\;\middle|\;
\begin{cases}
a_1a_2a_1=a_2a_1a_2,\,\,a_2a_3a_2=a_3a_2a_3,&\\
a_3a_4a_3=a_4a_3a_4,\,\,a_4a_1a_4=a_1a_4a_1,&\\
[a_1,a_3]=[a_2,a_4]=1&
\end{cases}
\right\rangle,
$$
which is isomorphic to the affine Artin group of type $\tilde{A}_3$.
\end{theorem}

\noindent\emph{Proof}. Write
$$F(x,y)=\frac{27}{4}y^4+2x(x^2-9)y^2-4(x^2-1)^2.$$
Since the leading coefficient of $F$ as a
polynomial in $y$ is a nonzero constant, the exceptional
values are precisely the zeros of its discriminant. A direct computation gives
$$\Delta_y F=-6912(x^2-1)^2(x^2+3)^6.$$
Thus $S=\left\{\pm 1,\pm\sqrt{-3}\right\}$, and the fiber of the projection $(x,y)\mapsto x$ over the base point $x_0=0$ is a four-punctured plane. Let $r=\left(\frac{16}{27}\right)^{1/4}$ and label the four punctures by
$y_1=-r$, $y_2=-\sqrt{-1}r$, $y_3=\sqrt{-1}r$ and $y_4=r$. For each $s\in S$, let $\gamma_s$ be a meridian based at $0$. The running tracks of these $y_i$'s along $\gamma_s$ will give the braid $f_*(\gamma_s)\in B_4$ that will be needed in the computation. To obtain this, note that $F(x,y)=0$ implies
$$y^2=\frac4{27}
\left(x(9-x^2)\pm(x^2+3)^{3/2}
\right).$$ This gives
\begin{equation}
\begin{aligned}f_*(\gamma_{1})&=\sigma_2^2\\
f_*(\gamma_{-1})&=\sigma_1\sigma_3\sigma_2^2\sigma_3^{-1}\sigma_1^{-1}\\
f_*\left(\gamma_{\sqrt{-3}}\right)&=(\sigma_1\sigma_3)^3\\
f_*\left(\gamma_{-\sqrt{-3}}\right)&=\sigma_2^{-1}(\sigma_1\sigma_3)^3\sigma_2.
\end{aligned}
\end{equation}
\begin{figure}
\centering
\begin{tikzpicture}[line cap=round,line join=round,>=Stealth]
\def\R{2.25}     
\def\k{0.28}     
\def\H{6.8}      
\def\a{0.92}     
\def\tau{0.18}   
\pgfmathsetmacro{\ry}{\k*\R}
\pgfmathsetmacro{\ay}{\k*\a}
\colorlet{colone}{blue!75!black}
\colorlet{coltwo}{red!75!black}
\colorlet{colthree}{teal!70!black}
\colorlet{colfour}{green!80!black}

\draw(\R,0) arc[start angle=0,end angle=180,x radius=\R,y radius=\ry];
\draw[dashed] (\R,-\H) arc[start angle=0,end angle=180,x radius=\R,y radius=\ry];

\draw (-\R,0) -- (-\R,-\H);
\draw (\R,0) -- (\R,-\H);

\draw[line width=1.05pt,colone] (-\R,0) -- (-\R,-\H);

\draw[line width=1.05pt,colfour] (\R,0) -- (\R,-\H);

\draw[dashed,line width=1.05pt,coltwo]
plot[domain=0:\tau,samples=40,variable=\t]
  ({0},{-\H*\t + \k*(\R - (\R-\a)*\t/\tau)});

\draw[dashed,line width=1.05pt,coltwo]
  plot[domain=\tau:{\tau + 0.25*(1-2*\tau)},samples=60,variable=\t]
  ({\a*sin(360*(\t-\tau)/(1-2*\tau))},
   {-\H*\t + \k*(\a*cos(360*(\t-\tau)/(1-2*\tau)))});

\draw[line width=1.05pt,coltwo]
  plot[domain={\tau + 0.25*(1-2*\tau)}:{\tau + 0.75*(1-2*\tau)},samples=120,variable=\t]
  ({\a*sin(360*(\t-\tau)/(1-2*\tau))},
   {-\H*\t + \k*(\a*cos(360*(\t-\tau)/(1-2*\tau)))});

\draw[dashed,line width=1.05pt,coltwo]
  plot[domain={\tau + 0.75*(1-2*\tau)}:{1-\tau},samples=60,variable=\t]
  ({\a*sin(360*(\t-\tau)/(1-2*\tau))},
   {-\H*\t + \k*(\a*cos(360*(\t-\tau)/(1-2*\tau)))});

\draw[dashed,line width=1.05pt,coltwo]
  plot[domain={1-\tau}:1,samples=40,variable=\t]
  ({0},{-\H*\t + \k*(\a + (\R-\a)*(\t-(1-\tau))/\tau)});
\draw[line width=1.05pt,colthree]
  plot[domain=0:\tau,samples=40,variable=\t]
  ({0},{-\H*\t - \k*(\R - (\R-\a)*\t/\tau)});
\draw[line width=1.05pt,colthree]
  plot[domain=\tau:{\tau + 0.25*(1-2*\tau)},samples=60,variable=\t]
  ({-\a*sin(360*(\t-\tau)/(1-2*\tau))},
   {-\H*\t - \k*(\a*cos(360*(\t-\tau)/(1-2*\tau)))});

\draw[dashed,line width=1.05pt,colthree]
  plot[domain={\tau + 0.25*(1-2*\tau)}:{\tau + 0.75*(1-2*\tau)},samples=120,variable=\t]
  ({-\a*sin(360*(\t-\tau)/(1-2*\tau))},
   {-\H*\t - \k*(\a*cos(360*(\t-\tau)/(1-2*\tau)))});

\draw[line width=1.05pt,colthree]
  plot[domain={\tau + 0.75*(1-2*\tau)}:{1-\tau},samples=60,variable=\t]
  ({-\a*sin(360*(\t-\tau)/(1-2*\tau))},
   {-\H*\t - \k*(\a*cos(360*(\t-\tau)/(1-2*\tau)))});

\draw[line width=1.05pt,colthree]
  plot[domain={1-\tau}:1,samples=40,variable=\t]
  ({0},{-\H*\t - \k*(\a + (\R-\a)*(\t-(1-\tau))/\tau)});

\draw (-\R,0) arc[start angle=180,end angle=360,x radius=\R,y radius=\ry];
\draw (-\R,-\H) arc[start angle=180,end angle=360,x radius=\R,y radius=\ry];
\fill[colone] (-\R,0) circle (1.7pt);
\node[text=colone,left=3pt] at (-\R,0) {$y_1=-r$};

\fill[coltwo] (0,\ry) circle (1.7pt);
\node[text=coltwo,above=2pt] at (0,\ry) {$y_2=ir$};

\fill[colthree] (0,-\ry) circle (1.7pt);
\node[text=colthree,below=2pt] at (0,-\ry) {$y_3=-ir$};

\fill[colfour] (\R,0) circle (1.7pt);
\node[text=colfour,right=3pt] at (\R,0) {$y_4=r$};
\fill[colone] (-\R,-\H) circle (1.7pt);
\node[text=colone,left=3pt] at (-\R,-\H) {$y_1$};

\fill[coltwo] (0,-\H+\ry) circle (1.7pt);
\node[text=coltwo,above=2pt] at (0,-\H+\ry) {$y_2$};

\fill[colthree] (0,-\H-\ry) circle (1.7pt);
\node[text=colthree,below=2pt] at (0,-\H-\ry) {$y_3$};
\fill[colfour] (\R,-\H) circle (1.7pt);
\node[text=colfour,right=3pt] at (\R,-\H) {$y_4$};
;
\end{tikzpicture}
\caption{running tracks along $\gamma_1$}
\label{figure 9}
\end{figure}
Indeed, along $\gamma_1$, the two imaginary roots $y_2,y_3$ coalesce at the node, with no intervening projection
crossings, see figure \ref{figure 9}. Similarly, along the path to $\gamma_{-1}$, the two real roots coalesce; the connecting path contributes the two negative crossings $\sigma_1^{-1}$ and $\sigma_3^{-1}$. At $\sqrt{-3}$, the pairs $(y_1,y_2)$ and $(y_3,y_4)$ coalesce at two cusps, with no projection crossings along the connecting path. Along the path to $-\sqrt{-3}$, one positive crossing $\sigma_2$ occurs before the two pairs coalesce. The local braids at a node and a cusp are respectively the square and the cube of a positive half-twist. Including the return paths gives the four braids above.

By Theorem \ref{theorem 2.2} we can compute that 
$$\begin{aligned}
[g_2,g_3]&=1=\left[g_2g_1g_2^{-1},g_4\right],\\
g_1g_2g_1&=g_2g_1g_2,\\
g_3g_4g_3&=g_4g_3g_4,\\
g_1g_3g_1&=g_3g_1g_3,\\
g_2g_4g_2&=g_4g_2g_4.
\end{aligned}$$
Define $a_3=g_2^{-1}g_4g_2$, $a_1=g_1$, $a_2=g_2$ and $a_4=g_3$, then the theorem follows directly. $\clubsuit$

For a regular value $\mathbf{c}\in\C^2\setminus\mathscr{C}$, by giving a partial compactification $\overline{\mu^{-1}(\mathbf{c})}$ of the generic fiber $\mu^{-1}(\mathbf{c})$, the periodic lattice $\Lambda_\mathbf{c}$ associated to $\overline{\mu^{-1}(\mathbf{c})}$ is also constructed in \cite{cushman2001complex}, it turns out to be a rank 3 lattice in $\C^2$. The complex Hamiltonian monodromy of $\Lambda$ around a specific loop in $\pi_1(\C^2\setminus\mathscr{C})$ was also computed in \cite[\S 6]{cushman2001complex}, it is conjugated to the matrix $M=\begin{pmatrix}1& -1 & 0\\ 0 & 1& 0 \\0& 0 & 1\end{pmatrix}$. Theorem \ref{theorem 4.1} shows that the complexified spherical pendulum has more fruitful information than the classical case, at least at the monodromy level. 

\renewcommand{\refname}{\textsc{References}}
\addcontentsline{toc}{section}{\textsc{References}}
\bibliographystyle{plain}
\footnotesize
\bibliography{references.bib}
\end{document}